\documentclass{article}
\usepackage[utf8]{inputenc}

\usepackage[utf8]{inputenc}
\usepackage[dvipsnames]{xcolor}
\usepackage{geometry}
\usepackage{xfrac}
\usepackage{amsfonts}
\usepackage{graphicx}
\usepackage{amsthm}
\usepackage{amsmath}
\usepackage{amssymb}
\usepackage{enumitem}
\usepackage{tikz-cd}
\usepackage{hyperref}
\usepackage{cleveref}
\usepackage{mathtools}
\usepackage[affil-it]{authblk}
\usepackage{ulem}

\graphicspath{ {./Figures/} }

\newcommand{\R}{\mathbb{R}}
\newcommand{\M}{\mathcal{M}}
\newcommand{\tnn}{\geq 0}
\newcommand{\Grkn}{Gr_{k,n}}
\newcommand{\nSk}{{{[n]}\choose{k}}}
\newcommand{\nCk}{{{n}\choose{k}}}
\newcommand{\tilingSet}{\mathbf{Til}}
\newcommand{\tilingSetN}{\mathbf{Til_n}}
\newcommand{\strandDiagSet}{\mathbf{Diag}}
\newcommand{\strandDiagSetN}{\mathbf{Diag_n}}
\newcommand\restrict[1]{\raisebox{-.5ex}{$|$}_{#1}}

\newtheorem{theorem}{Theorem}[section]
\newtheorem{proposition}[theorem]{Proposition}
\theoremstyle{definition}
\newtheorem{definition}[theorem]{Definition}
\newtheorem{remark}[theorem]{Remark}
\newtheorem{lemma}[theorem]{Lemma}
\newtheorem{corollary}[theorem]{Corollary}
\newtheorem{example}[theorem]{Example}

\newtheorem{preTheorem}{Theorem}

\newenvironment{acknowledgements} {\begin{abstract}}{\end{abstract}}

\title{Parametrizing positroid cells using bicolored tilings}
\author{Joel Costa}
\affil{School of Mathematics, University of Leeds, United Kingdom}

\date{March 14, 2022}

\begin{document}

\maketitle

\begin{abstract}
    Bicolored tilings are given by a collection of smooth curves in a disk with a coloring map on the tiles these curves form. Postnikov diagrams can be viewed as the image of certain bicolored tilings under the Scott map. We introduce a reduction technique on bicolored tilings, and show that a tiling maps to a Postnikov diagram if and only if it is reduced. We then use bicolored tilings to parametrise positroid cells in the Grassmannian, and use the reduction, along with another transform, to generate tilings associated to lower-dimensional positroids cells. We also show that the parametrisation of such a cell can be derived from the parametrisation of the higher-dimensional cell.
\end{abstract}

\section{Introduction}

The Grassmannian $\Grkn$ is the space of $k$-dimensional subspaces in an $n$-dimensional vector space $V$. For our intents, that space is $V = \R^n$. We may represent a point $W \in \Grkn$ as a $k \times n$ matrix with entries in $\R$, with $W$ being the row space of that matrix. The totally non-negative Grassmannian $\Grkn^{\tnn}$ is the part of $\Grkn$ consisting of matrices for which all maximal minors are non-negative. $\Grkn^{\tnn}$ can be stratified into positroid cells \cite{Postnikov}, which we parametrise using bicolored tilings. Our main results are:

\begin{preTheorem}
    Reduced (bicolored) tilings of type $(k,n)$ up to tiling equivalence are in bijection with positroid cells of the totally non-negative Grassmannian $\Grkn^{\tnn}$.
\end{preTheorem}

\begin{preTheorem}
    Let $T$ be a tiling of type $(k,n)$ with permutation $\pi$, and let $\alpha \in A$ be an angle in $T$. Let $T' \vcentcolon = d_{\alpha}(T)$ be the degeneration of $T$ at $\alpha$, and let $P=P(\beta)_{\beta \in A}$ be the parametrisation of $T$. Then
    \begin{enumerate}[label = $\bullet$]
        \item $T'$ is of type $(k,n)$.
        \item $T'$ has decorated permutation $\pi' = (i \,\, j) \circ \pi$ for some $i,j \in [n]$.
        \item $T'$ parametrises the positroid cell $S_{\pi'}$ by $P \restrict{\alpha = 0}$.
        \item $T < T'$ and $\dim S_T < \dim S_{T'}$.
    \end{enumerate}
\end{preTheorem}

\begin{acknowledgements}
    The author would like to thank Karin Baur her help and insightful comments during the writing of this paper and the work that preceded it. The author was supported by Royal Society Wolfson Fellowship 180004.
    
    The author would also like to thank the Isaac Newton Institute for Mathematical Sciences, Cambridge, for support and hospitality during the programme \textit{Cluster algebras and representation theory} where work on this paper was undertaken. This work was supported by EPSRC grant no EP/K032208/1.
\end{acknowledgements}

\section{Background}

We start by recalling (bicolored) tilings, Postnikov diagrams, as well as the Scott map which links the two. More details can be found in \cite{Costa} where bicolored tilings were first introduced.

First we generalise edges to edges between any finite number of vertices, which can be viewed as polygons (up to isotopy) with those vertices as endpoints.

\begin{definition} \cite{Costa} \label{edgeDefinition}
    Let $S$ be a $2$-dimensional connected oriented surface with boundary, with $n$ distinct boundary vertices, enumerated $\{1,\dots,n\}$, and $x_1,\dots,x_m$ internal vertices, for $m \geq 0$. Let $V$ be the set of vertices, both boundary and internal. An \textit{edge} $e=(v_1,\dots,v_r)$, $r > 0$, is a finite sequence of vertices $v_i \in V$ such that
    \begin{enumerate}[label = $(\roman*)$]
        \item There are no repetitions of vertices in $v_1,\dots,v_r$, with the exception of boundary vertices which may appear exactly twice in consecutive order (with the convention that $v_r$ and $v_1$ are consecutive).
        \item There is a collection of smooth curves $\gamma_1, \dots, \gamma_m$ on the surface such that $\gamma_i$ has endpoints $v_i$ and $v_{i+1}$ (with the convention that $v_{r+1} = v_1$), and such that no two curves intersect, other than at the endpoints of any two consecutive curves $\gamma_i$ and $\gamma_{i+1}$.
        \item There are no vertices of $V$ in the interior of the disk with boundary $\bigcup \gamma_i$.
     \end{enumerate}
     We call $\partial e = \bigcup \gamma_i$ the \textit{boundary} of $e$, and $\partial^2 e = \{v \mid v \in e\}$ the \textit{endpoints} of $e$. Edges are treated up to cyclical shift of the sequence.
\end{definition}

\begin{remark}
    We draw an edge with boundary curves $\gamma_1,\dots,\gamma_r$ by shading the area between those curves. We also refer to such an edge as a \textit{black tile}. Edges between $2$ distinct vertices can be drawn as either an arc between the vertices or black digons with the vertices as endpoints (and are considered equivalent), and are also called \textit{simple edges}. An edge consisting of a single vertex is drawn as a black $1$-gon (a loop whose interior is shaded).
\end{remark}

\begin{example}
    We consider the disk $D_6$ with $6$ boundary vertices and $3$ internal vertices $x_1,x_2,x_3$.
    \begin{center}
        \includegraphics[scale=0.45]{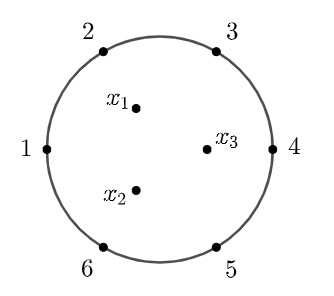}
    \end{center}
    Then for $b_i = (i,i+1)$ (for $i=1,\dots,6$), the following are examples of edges.
    \begin{enumerate} [label=$\cdot$]
        \item $b_1,\dots,b_6$, $e_1 = (2,x_1)$, $e_2 = (x_1,x_2,x_3)$, $e_3 = (5,x_3)$
        \item $b_1,b_2,b_3,b_5,b_6$, $e_4 = (2,x_1)$, $e_5=(3,4,5,x_1)$
        \item $b_2,\dots,b_5$, $e_6=(1,2,x_1)$, $e_7=(1,x_1,6), e_8=(x_2)$
    \end{enumerate}
    \begin{center}
        \includegraphics[scale=0.4]{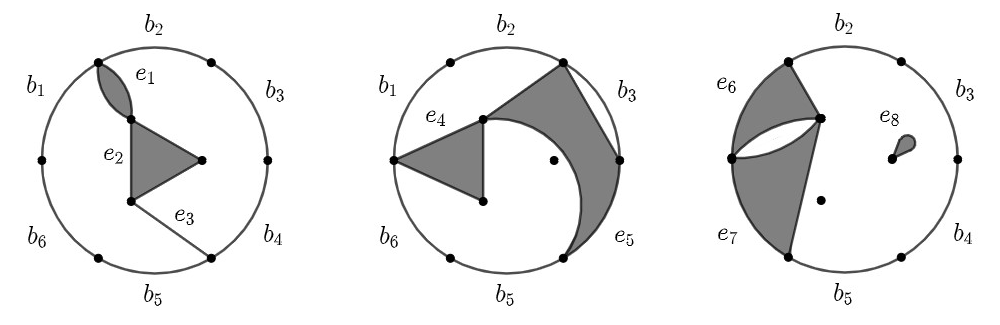}
    \end{center}
    We note that if a boundary vertex $i$ appears twice in the sequence that describes an edge, the curve that has both endpoints $i$ forms a white loop at the boundary, as seen in the second example. If two vertices $u,v$ appear as consecutive vertices in two (or more) edges, this creates (multiple) white digons, as seen in the second and third example.
\end{example}

\begin{definition}
    Let $S$ be a $2$-dimensional connected oriented surface with boundary, with $n$ distinct boundary vertices, enumerated $\{1,\dots,n\}$, and $x_1,\dots,x_m$ internal vertices. We denote the set of these vertices by $V$. A \textit{bicolored} tiling $T = (S,V,E)$ is the surface $S$ equipped with a finite collection of edges $E$ such that
    \begin{enumerate} [label = $\cdot$]
        \item the black tiles representing any two distinct edges only intersect at their common endpoints.
        \item for any boundary segment between two consecutive boundary vertices $u,v$, there is exactly one black tile intersecting that whole boundary segment, including the endpoints $u,v$.
    \end{enumerate}
    We define tilings up to equivalence given by the following two local transformations
    \begin{enumerate}[label = $\cdot$]
        \item Hourglass equivalence:
            \begin{center}
                \includegraphics[scale=0.65]{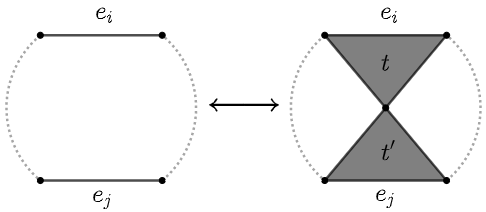}
            \end{center}
        \item Digon equivalence:
            \begin{center}
                \includegraphics[scale=0.65]{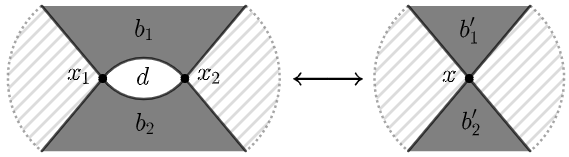}
            \end{center}
            where $x_1$ and $x_2$ are not simultaneously boundary vertices.
    \end{enumerate}
    The set of tilings up to tiling equivalence, i.e. hourglass/digon equivalence, is denoted $\tilingSet$. The set of tilings with $n$ boundary vertices is denoted $\tilingSet_n$.
    Furthermore, we define the \textit{mutation/flip} of a simple edge inside a quadrilateral within the tiling, as described in the following figure
    \begin{center}
        \includegraphics[scale=0.60]{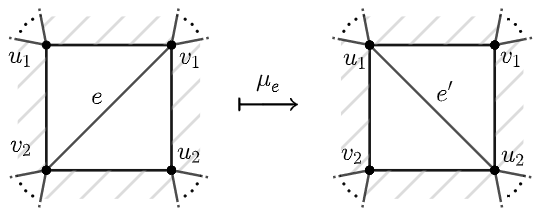}
    \end{center}
    We say that two tilings are flip/mutation-equivalent if one can be obtained from the other by a finite sequence of flips and tiling equivalences. We denote $\overline{\tilingSet}$ the set of flip/mutation-equivalence classes of tilings in $\tilingSet$. For simplicity, we will identify $T$ with its flip-equivalence class $\overline{T}$ for the rest of this paper, unless otherwise specified.
\end{definition}

\begin{definition}
    Let $T$ be a tiling of a surface $S$, let $t_e$ be the black tiles corresponding to $e$, for any edge $e \in E$. The \textit{faces} $T$ are the connected components of $S \setminus (\bigcup t_e)$. We also call the faces of a tiling \textit{white tiles}. The set of faces of $T$ is denoted $F$.
    
    An \textit{angle} $\alpha$ of $T$ is a quadruple $(v,e_1,e_2,f) \in V \times E \times E \times F$, such that $e_1$ and $e_2$ intersect. Informally, we choose a vertex around which the angle $\alpha$ lies, and a face inside which $\alpha$ lies. However, this is not enough in some cases, such as when there is a $1$-edge, and thus we need to specify between which two edges $\alpha$ lies. The set of angles of $T$ is denoted $A$. 
\end{definition}

\begin{definition} \label{unitedTilings}
    Let $\mathcal{C} = S_{1} \cup \dots \cup S_{r}$ be a finite union of $2$-dimensional connected oriented surfaces with $n_1,\dots,n_r$ boundary vertices, such that the intersection of any two surfaces $S_{i}$ and $S_{j}$ is either empty or a common boundary vertex. Then a \textit{tiling} $\mathcal{T}$ of $\mathcal{C}$ is the union of tilings $T_1 \cup \dots \cup T_r$ such that $T_i=(S_i,V_i,E_i)$ is a tiling of $S_{i}$, where $T_i \cup T_j = (S_i \cup S_j, V_i \cup V_j, E_i \cup E_j)$.
\end{definition}

\begin{example}
    The following is a tiling of $D_6$, a tiling that is equivalent to it, and a tiling that is flip-equivalent to it.
    \begin{center}
        \includegraphics[scale=0.4]{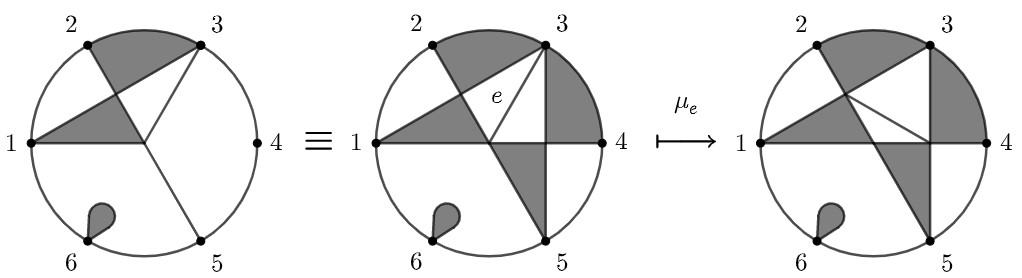}
    \end{center}
    We will usually draw $D_n$ as an $n-gon$, as we work with tilings up to homotopy. This would make our tiling of $D_6$ look as follows
    \begin{center}
        \includegraphics[scale=0.4]{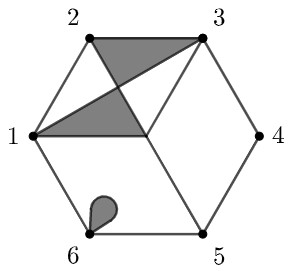}
    \end{center}
\end{example}

\begin{remark}
    The number of internal vertices varies between equivalent tilings. The hourglass equivalence adds or removes the middle vertex when we go right or left in the above depicted transformation, respectively. The digon equivalence contracts two vertices into one from left to right in the depiction above.
\end{remark}

\begin{definition} \label{subtilingDefinition}
    Let $T$ be a tiling of a surface $S$ with vertex set $V$ and edge set $E$. Let $v_1,\dots,v_m$ be vertices of $T$ such that for any $i \in [m]$, there is a curve $\gamma_i$ with endpoints $v_i$, $v_{i+1}$ (with $v_{m+1} =v_1$), as described in \Cref{edgeDefinition}. Let $S' \subset S$ be a surface with boundary $\bigcup \gamma_i$ and boundary vertices $v_1,\dots,v_m$. Then we define a tiling $T'$ of $S'$ as follows
    \begin{enumerate}[label = $\cdot$]
        \item edge set $E' = \{e \cap S' \,|\, e \in E \text{ and } e \cap S' \neq \emptyset \}$,
        \item vertex set $V' = V \cap S'$,
        \item boundary vertices $\partial V' = \{v_1,\dots,v_m\} \cup (\partial V \cap \partial S)$.
    \end{enumerate}
    We say that $T'$ is a \textit{subtiling} of $T$ under $S'$. Let $\tilde S = S \setminus \text{int}(S') \subset S$ be the surface with boundary $\partial \tilde{S} = \overline{\partial S \cup \partial S' \setminus (\partial S \cap \partial S')}$. Then we call the subtiling $\tilde{T}$ under $\tilde{S}$ the \textit{remainder} of $T$ under $S'$. We also denote $\tilde{T} = T \setminus T'$.
\end{definition}

For the rest of this paper, unless otherwise specified, we only consider tilings for the disk $S=D_n$ with $n$ boundary vertices, and the boundary vertices are labeled $1,\dots,n$ in clockwise order.

Next, we want to define \textit{Postnikov diagrams}. Before we give the definition in \cite{Postnikov}, we define a more general notion which allows for so-called \textit{"bad double crossing"}.

\begin{definition}
    Consider a disk with $n$ vertices drawn on its boundary, labeled by the elements in $\{1,\dots,n\}$, in clockwise order. An \textit{(alternating) curve diagram} consists of a finite collection of oriented curves, such that each curve is either a closed cycle or has boundary vertices as endpoints, in which case we call it a strand, with every boundary vertex having exactly one incoming and outgoing strand, and satisfying the following conditions:
    \begin{enumerate}[label=(\roman*)]
        \item
            A curve does not cross itself in the interior of the disk.
        \item
            No three curve cross in one single point.
        \item
            All crossings are transversal (left figure as opposed to right figure).
            \begin{center}
                \includegraphics[scale=0.50]{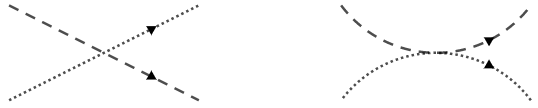}
            \end{center}
        \item
            There are finitely many crossings between curves.
        \item
            Following any curve in one direction, the curves that intersect it must alternate in orientation.
            \begin{center}
                \includegraphics[scale=0.5]{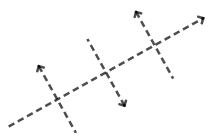}
            \end{center}
    \end{enumerate}
    We define alternating curve diagrams up to equivalence of two local transformations, namely twisting and untwisting oriented lenses inside the disk or on the boundary
    \begin{center}
        \includegraphics[scale=1.5]{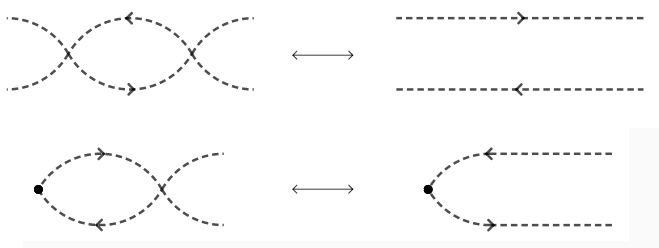}
    \end{center}
    We call the transformations from left to right a \textit{reduction}, and diagrams to which no further reduction can be applied \textit{reduced diagrams}. If two curve diagrams $D_1,D_2$ are equivalent up to these transformations, we write $D_1 \equiv D_2$. The set of alternating curve diagrams in a disk with $n$ boundary vertices up to equivalence is denoted $\strandDiagSetN$. The set of all alternating curve diagrams is denoted $\strandDiagSet = \bigcup \strandDiagSetN$.
    Furthermore, we treat every diagram up to isotopy with the boundary vertices fixed. When necessary, we denote this equivalence $\sim$.
    For any $i \in \{1,\dots,n\}$, the strand that starts at the boundary vertex $i$ is denoted $\gamma_i$. Any curve diagram has exactly $n$ strands $\gamma_1,\dots,\gamma_n$, and may have closed cycles in the interior of the disk.
\end{definition}

\begin{definition} \cite{Postnikov}
    A \textit{Postnikov diagram}, or (alternating) strand diagram is an alternating curve diagram such that
    \begin{enumerate}[label=(\roman*)]
        \item
            There are no closed cycles, i.e. every curve is a strand. Equivalently, there are as many curves as there are boundary vertices.
        \item
            If two strands cross at two points $A$ and $B$, then one strand is oriented from $A$ to $B$, and the other from $B$ to $A$ (left figure as opposed to right figure). In other words, no two strands create unoriented lenses, or more informally,there are no \textit{bad double crossings}.
            \begin{center}
                \includegraphics[scale=0.50]{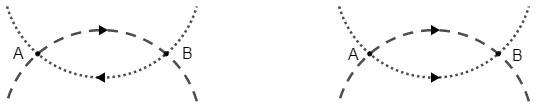}
            \end{center}
    \end{enumerate}
\end{definition}

\begin{definition} \cite{Postnikov}
    A \textit{decorated permutation} $\overline{\pi}$ of $\{1,\dots,n\}$ is a pair $(\pi,c)$ consisting of a permutation $\pi$ of $\{1,\dots,n\}$ and a map $c$ that maps any fixed point of $\pi$ to an element in $\{-1,1\}$.
\end{definition}

The function $c$ is a colouring of the fixed points of $\pi$. Any strand diagram defines a permutation $\pi$ of $\{1,\dots,n\}$ where $\pi(i) = j$ when $\gamma_i$ ends at the boundary vertex $j$. For any fixed point $i$ where $\gamma_i$ is oriented clockwise, $col(i) = 1$. Otherwise $col(i)=-1$. We call $i$ the \textit{source (vertex)} and $\pi(i)=j$ the \textit{target (vertex)} of $\gamma_i$.

\begin{definition}
    Let $\Gamma$ be a curve diagram with $n$ boundary vertices. Let $\pi$ be the a decorated permutation of $\Gamma$, written as $\pi: i \longmapsto i + s(i)$, where $s(i) \in \{0, \dots, n\}$, such that $s(i) = n$ exactly then when $col(i)=1$, i.e. $\gamma_i$ is a clockwise loop starting and ending at $i$. Then the \textit{rank} of $\Gamma$ is given by
    \[
        k = \frac{1}{n} \sum_{i=1}^{n} s(i)
    \]
    We call $(k,n)$ the \textit{type} of $\Gamma$.
\end{definition}

\begin{definition} \cite{Scott} \label{ScottMap}
    We define the \textit{Scott map}
    $$S: \tilingSetN \longrightarrow \strandDiagSetN, T \longmapsto D$$
    to be the map such that
    \begin{enumerate}[label=$\cdot$]
        \item any white tile is mapped to a configuration consisting of $m$ curve segments, where $m$ is the size of the tile, following around the border in a counter-clockwise orientation. For example:
        \begin{center}
            \includegraphics[scale=0.54]{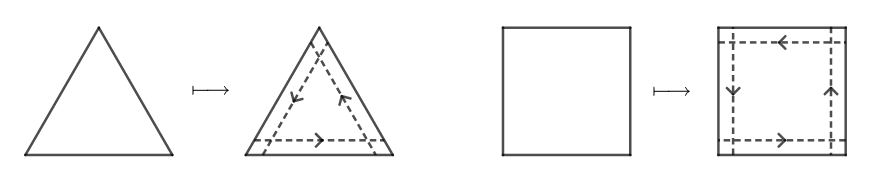}
        \end{center}
        \item any black tile is mapped to a configuration consisting of $m$ curve segments, where $m$ is the number of vertices of the tile, such that each curve forms an arc around a vertex inside the tile in a clockwise orientation. For example:
        \begin{center}
            \includegraphics[scale=0.54]{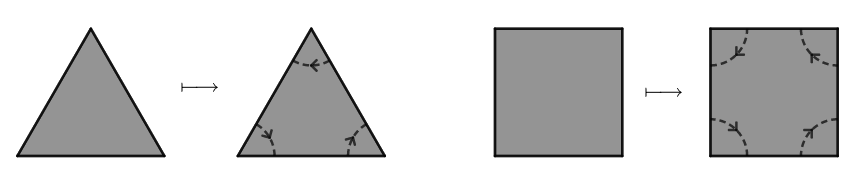}
        \end{center}
        \item If two tiles are adjacent, join the pairs of oriented curves segments along their shared boundary. The oriented curves obtained from concatenating all curve segments make up the full curves of the diagram. One can check that these are consistently oriented. Indeed, the only intersections of curves occur in white tiles. Following a single curve, these intersections always come in pairs of two, the first being from left to right, and the second being from right to left. Thus the intersections keep alternating as the curve passes through white tiles.
        \item The curves join at the boundary, i.e. for any boundary vertex, we take the two curves that intersect the boundary on either side of the vertex closest to it and join them.
    \end{enumerate}
\end{definition}

\begin{definition}
    The \textit{rank} of a tiling $T$ is defined as the rank of the diagram $S(T)$. Similarly, the type of $T$ is defined as the type of $S(T)$. 
\end{definition}

We are mainly interested in tilings whose image under the Scott map are Postnikov diagrams. However, not all tilings map to a Postnikov diagram, as the following example illustrates.
\begin{center}
    \includegraphics[scale=1.85]{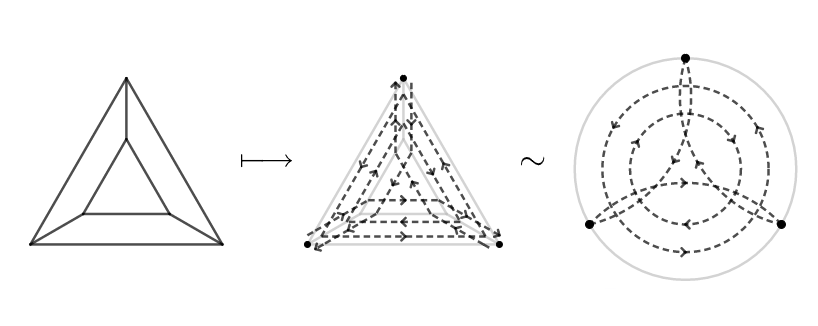}
\end{center}

In order to obtain Postnikov diagrams, we will later define a transformation on tilings that reduces the number of edges (and black tiles in general), allowing us to produce a tiling that maps to a Postnikov diagram, and associate the original tiling to that diagram (see \Cref{tilingReduction}). This reduction is a reasonable choice as it preserves combinatorial properties related to the positroid cell associated to the decorated permutation of the resulting Postnikov diagram.

Before we move on to the next section where we define this reduction, we recall one more notion from \cite{Postnikov}.

\begin{definition} \cite{Postnikov}
    A \textit{plabic graph} $G$ is a planar bipartite graph inside a disk $D$ with $n$ designated vertices of degree $1$ on the boundary of $D$. The internal vertices in either part of the graph are colored black or white, respectively.
\end{definition}

\begin{remark} \cite{BaurMartin} \label{stellarReplacementMap}
    Tilings naturally map to plabic graphs by the map $\Phi$ as follows

    \begin{enumerate}[label = $\cdot$]
        \item Any vertex $v$ of $T$ is mapped to a white vertex $\Phi(v)$ of $G$.
        \item Any face $f$ of $T$ is mapped to a black vertex $\Phi(f)$ of $G$.
        \item If $v \in \partial f$, then there is an edge between $\Phi(v)$ and $\Phi(f)$ in $G$. In other words, any angle $\alpha \in A$ around $v$ and in $f$ is mapped to an edge $\Phi(\alpha)$ in $G$ between $\Phi(v)$ and $\Phi(f)$.
        \item We draw a circle with $n$ vertices labeled $1'$ to $n'$ in clockwise order. These will be the boundary vertices of $G$. Then for any boundary vertex $i$ in $T$, we add an edge from $\Phi(i)$ and $i'$ in $G$.
    \end{enumerate}
    
    Then the resulting construction is the plabic graph $\Phi(T)$.
    \begin{center}
        \includegraphics[scale=0.4]{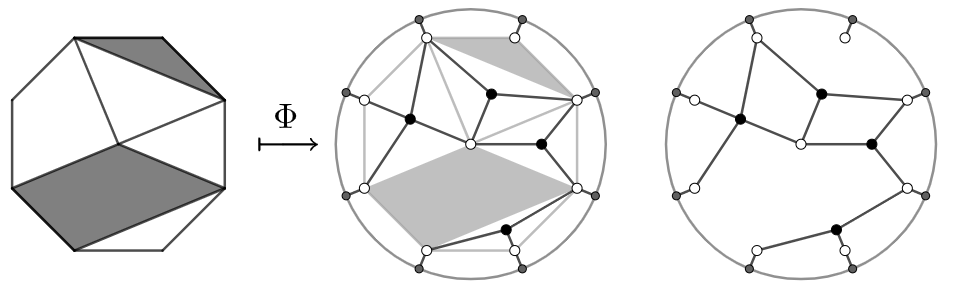}
    \end{center}
    We call $\Phi$ the \textit{stellar-replacement map}. It was introduced in \cite{BaurMartin} for the case of unicolored tilings (bicolored tilings whose edges are all simple). The map naturally extends to the bicolored set-up.
\end{remark}

\section{Reduction of bicolored tilings}

We define another transformation on tilings called a reduction. This transformation is linked to parallel edge reductions of plabic graphs in \cite{Postnikov}[12.4, p.43] and allows us to describe tilings that map to Postnikov diagrams. Furthermore, reductions on tilings preserve the positroid cell associated to tilings, as we will see in Section 4.

\begin{definition} \label{tilingReduction}
    Let $T$ be a tiling, and let $e \in E$ be a black $1$-gon whose only neighboring tile is a white tile. Then the tiling $T' = T-e$ is called a \textit{reduction} of $T$. We denote this reduction $R_e$.
    \begin{center}
        \includegraphics[scale=0.6]{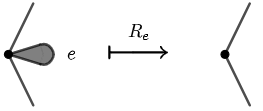}
    \end{center}
\end{definition}

\begin{definition}
    A tiling $T$ is said to be \textit{reduced} if no reduction $R_e$ can be applied to any tiling in its mutation-equivalence class.
\end{definition}

It can be hard to see when a tiling is reduced. However, we will see in \Cref{reducedTilingPostnikov} that reduced tilings map to Postnikov diagrams. Thus, if a tiling does not generate a Postnikov diagram under the Scott map, we know it is not reduced, and try to find a tiling in its mutation-equivalence class to reduce it. We do this until we obtain a reduced tiling that maps to a Postnikov diagram. Later, we will see that the reduced tiling preserves some combinatorial properties of the initial tiling (particularly in \Cref{invariancesParametrisation}).

\begin{proposition}
    Let $T$ be a tiling of permutation $\pi$, $e$ be a black $1$-gon, and $T' = R_e(T)$. Let $\gamma_i$ and $\gamma_j$ be the strands in $T$ that pass through $e$ and go around $e$, respectively. Then the permutation $\pi'$ of $T'$ is given by 
    \[
    \pi'(l) =
    \begin{cases}
        \pi(l), & \text{if } l \neq i,j\\
        \pi(j) & \text{if } l = i \\
        \pi(i) & \text{if } l = j
    \end{cases}
    \]
    In other words, $\pi' = (\pi(i) \,\, \pi(j)) \pi$.
    \begin{proof}
        We observe the effect that reductions have from the diagram $S(T)$ to $S(T')$.
        \begin{center}
            \includegraphics[scale=0.65]{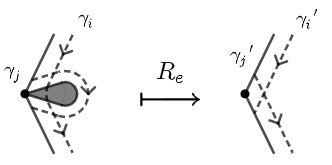}
        \end{center}
        We observe that the strands $\gamma_i$ and $\gamma_j$ simply swap target from $S(T)$ to $S(T')$.
    \end{proof}
\end{proposition}

\begin{example}
    The following tiling of a hexagon is equivalent (by hourglass equivalence) to a tiling with a black $1$-gon $e$.
    \begin{center}
        \includegraphics[scale=0.5]{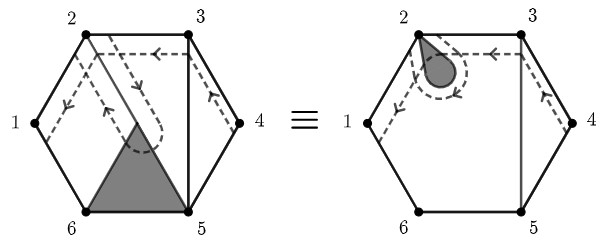}
    \end{center}
    If $\pi$ is the permutatation of these tilings, we can see that $\pi(2) = 2$ and $\pi(4) = 1$. By applying $R_e$ to the tiling, we obtain a tiling
    \begin{center}
        \hspace{0.05cm} \includegraphics[scale=0.5]{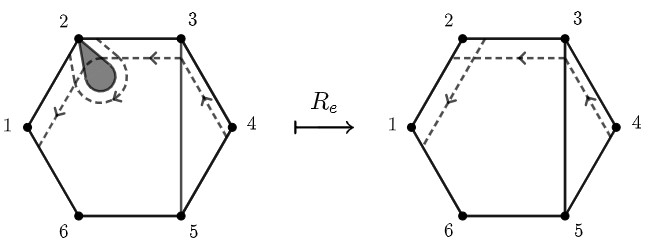}
    \end{center}
    If $\pi'$ is the permutation of the resulting tiling, then $\pi'(2) = 1$ and $\pi'(3) = 2$. It is easy to verify that for any $i\neq 2,3$, $\pi(i)=\pi'(i)$.
\end{example}

\begin{proposition}
    Let $T$ be a tiling, $e$ be a black $1$-gon, and $T' = R_e(T)$. Then $T$ and $T'$ have the same type.
    \begin{proof}
        $T$ and $T'$ have the same number of boundary vertices, so it remains to show that $T$ and $T'$ have the same rank. Let $\gamma_i$ and $\gamma_j$ be the strands in $T$ that pass through $e$ and go around $e$, respectively (as shown in the previous figure). Let $\gamma_i'$ and $\gamma_j'$ be the strands that intersect at the angle in $T'$ where $e$ was. If $\pi$ is the permutation of $T$, then the permutation $\pi'$ of $T'$ is given by
        \[
        \pi'(l) =
        \begin{cases}
            \pi(l), & \text{if } l \neq i,j\\
            \pi(j) & \text{if } l = i \\
            \pi(i) & \text{if } l = j
        \end{cases}
        \]
        Then, $s'(i) = s(i) + \pi(j) - \pi(i)$, and $s'(j) = s(j) + \pi(i) - \pi(j)$
        Then the rank $k'$ of $T'$ is given by
        \[
            k' = \frac{1}{n} \sum_{l=1}^n s'(l) = \frac{1}{n} \sum_{l=1}^n s(l) + [\pi(j) - \pi(i) + \pi(i) - \pi(j)] = k
        \]
    \end{proof}
\end{proposition}

The following proposition generalises the result in \cite[3.5]{BaurMartin} which states that - in the classic setup for tilings - tilings map to Postnikov diagrams under the Scott map. Once we introduce internal vertices and black tiles, this statement does not hold for every tiling anymore. An equivalent statement for plabic graphs is given in \cite[14.2]{Postnikov}.

\begin{proposition} \label{reducedTilingPostnikov}
    $T$ is a reduced tiling if and only if $S(T)$ is a Postnikov diagram.
    \begin{proof}
        Let $T$ be non-reduced. Then there is $T' \equiv T$ such that $T'$ has a $1$-edge $e$. Then, locally, $e$ maps to a bad double crossing under the Scott map, thus $S(T)$ is not a Postnikov diagram. Hence, if $S(T)$ is a Postnikov diagram, then $T$ is a reduced tiling.
        
        Let $T$ be a reduced tiling, and $G=\Phi(T)$ the corresponding plabic graph as described in \Cref{stellarReplacementMap}. If $\mathcal{F}$ is the set of faces of $G$, then $|E| = |\mathcal{F}|$. Assume $S(T)$ is not a Postnikov diagram, then by \cite[14.12]{Postnikov}, $G$ is not reduced. Then by \cite[12.5]{Postnikov}, $|\mathcal{F}|$ in the movement-reduction class of $G$. Then $|E|$ is not minimal within the reduction-equivalence class of $T$, and thus $T$ is not reduced, which is a contradiction. Hence, $S(T)$ is a Postnikov diagram.
    \end{proof}
\end{proposition}

\section{Parametrising positroid cells}

We will recall the definion of the totally non-negative Grassmannian and its decomposition into positroid cells. We then introduce a way to parametrise these cells using tilings. We show that each tiling in $\tilingSet$ gives us a different positroid cell.

\begin{definition}
    The \textit{Grassmannian} $\Grkn$ of \textit{type} $(k,n)$ is the set of $k$-dimensional subspaces in an $n$-dimensional vector space $\mathbb{V}$. Here, $\mathbb{V} = \R^n$.
    
    A point $V \in \Grkn$ can be described by a full-rank $k \times n$-matrix $M$, with $V$ being the row-space of $M$. The row-space of $M$ is invariant under left action by a non-singular $k \times k$-matrix. Thus, we can identify the Grassmannian as
    \[
        \Grkn = GL_k \backslash Mat_{k \times n}
    \]
    where $Mat_{k \times n}$ is the set of full-rank $k \times n$-matrices.
\end{definition}

We can embed $\Grkn$ into the projective space $\mathbb{P}^{\nCk-1}$ by setting a coordinate for any $k$-subset $I$ of $[n] := \{1,\dots,n\}$
\[
    \Delta_I = \Delta_I(M)
\]
where $\Delta_I(M)$ is the minor of the matrix composed of the column vectors of $M$ enumerated by $I$. Then the collection $(\Delta_I)_{I \in \nSk}$ gives us projective coordinates for $V$.

The \textit{totally non-negative Grassmannian} $\Grkn^{\tnn}$ is the subset of subspaces in $\Grkn$ for which all the projective coordinates are all non-negative up to simultaneous scaling with a factor $\lambda \neq 0$.

$\Grkn^{\tnn}$ can be decomposed into \textit{positroid cells} $S_{\pi}$ (described in \cite{Postnikov}) indexed by decorated permutations
\[
    \Grkn^{\tnn} = \bigsqcup_{\pi \in S_n^{*}} S_{\pi}
\]
where $S_n^{*}$ is the set of decorated permutations of $[n]$. For this paper, a cell $S_{\pi}$ is best thought of as the set of all points in $\Grkn$ whose maximal minors $\Delta_I$ are exactly non-zero for the same $k$-subsets $I$. That is, if $M,N \in S_{\pi}$, then $\Delta_I M = 0 \Leftrightarrow \Delta_I N = 0$. The ability to index this decomposition by decorated permutations is given by \cite[Theorem 17.1]{Postnikov}.

\vspace{0.6cm}

For the rest of this paper, unless otherwise specified, the set of vertices, edges, faces, and angles of $T$ will be denoted $V$, $E$, $F$, and $A$, respectively. If we have a second tiling $T'$, we use $V'$,$E'$,$F'$, and $A'$, respectively. We denote the angle near $v \in V$ inside the face $f$ by $\alpha_f^{v}$. In that case, we also denote $v = v(\alpha)$ and $f = f(\alpha)$.

\begin{definition} \label{matchingDefinition}
    A \textit{matching} $m \subset A$ of a tiling $T$ is a choice of angles of $T$ such that.
    \begin{enumerate}[label = (\roman*)]
        \item Each face is matched exactly once, i.e. for any two angles $\alpha, \beta \in m: f(\alpha) \neq f(\beta)$, and for any face $f$ of $T$, there is $\alpha \in m$ such that $f=f(\alpha)$.
        \item Each vertex is matched at most once, i.e. for any two $\alpha, \beta \in m: v(\alpha) \neq v(\beta)$.
        \item Each internal vertex is matched exactly once, i.e. on top of the second condition, for any internal vertex $v$ of $T$, there is $\alpha \in m$ such that $v=v(\alpha)$.
    \end{enumerate}
    We denote $\partial m = \{v \in V \, | \, v \text{ is a boundary vertex and not matched in } m\}$ the \textit{boundary} of $m$. The set of matchings of a tiling $T$ is denoted $\M (T)$.
\end{definition}

\begin{example} \label{runningExampleHex}
    Consider the following tiling $T$ of type $(4,6)$ with one internal vertex and $12$ angles.
    \begin{center}
        \includegraphics[scale=0.4]{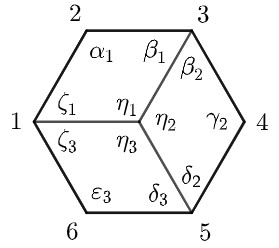}
    \end{center}
    Then three examples of matchings of $T$ would be
    \begin{center}
        \includegraphics[scale=0.4]{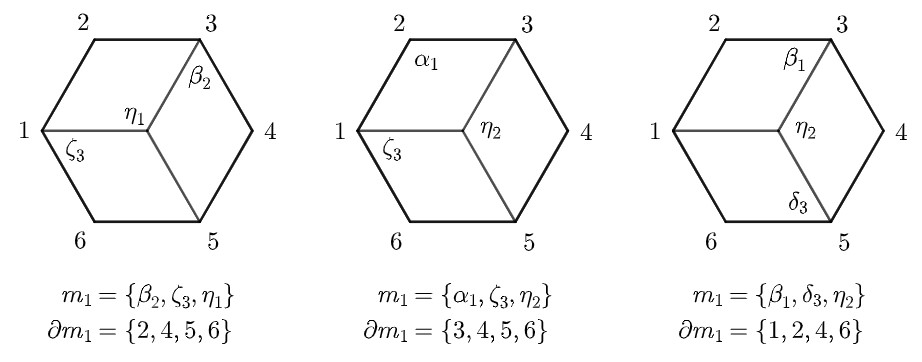}
    \end{center}
    For simplification, we often denote the $k$-subsets $I$ simply by concatenating their elements, that is $\partial m_1 = 2456$, $\partial m_2 = 3456$, and $\partial m_3 = 1246$.
\end{example}

\begin{remark}
    The definition of matchings of bicolored tilings maps onto the definition of almost perfect matching of plabic graphs in \cite{Williams} by the stellar-replacement map $\Phi$ described in \Cref{stellarReplacementMap}. The angles of a tiling map to edges in the corresponding plabic graph. This also gives us the following proposition.
\end{remark}

\begin{proposition}
    If $T$ is a tiling of type $(k,n)$ with a matching $m$, then $|\partial m| = k$.
\end{proposition}

Using this proposition, we can deduce the following.

\begin{proposition}
    Let $T$ be a tiling with a matching of type $(k,n)$. Then $k = |V|-|F|$.
    \begin{proof}
        If $m$ is a matching of $T$, then every face is matched to exactly one vertex. Then, $k = |\partial m|$ is given by the number of boundary vertices that are not matched, which is the number of total vertices minus the number of faces of the tiling.
    \end{proof}
\end{proposition}

\begin{definition} \label{tilingParametrisation}
    We consider any matching $m$ of $T$ as a monomial given by the product of its elements. Then for any $k$-subset $I$ of $[n]$, we set
    \[
        \Delta_I = \sum_{\substack{m \in \M (T)\\ \partial m = I}} m
    \]
    where $M(T)$ is the set of all matchings of $T$. The \textit{positroid cell} $S_T$ \textit{associated to} $T$ is given by all the points $(\Delta_I)_{I \in \nSk}$ for which the parameters $\alpha \in A$ are all positive, i.e.
    \[
        S_T = \{ (\Delta_I)_I \, | \, \alpha \in \R_{>0} \,\,\, \forall \alpha \in A \}
    \]
    We call $P_T = (\Delta_I)_I$ the parametrisation, and $\Delta_I$ the \textit{Plücker coordinates} of $T$ and $S_T$. The closure $\overline{S_T}$ of a positroid cell $S_T$ is given by 
    \[
        \overline{S_T} = \{(\Delta_I) \mid \alpha \in \R_{\geq 0} \,\,\, \forall \alpha \in A \}
    \]
    The closure of one cell $S_T$ is nested in another cell $S_{T'}$ if the zero coordinates of $S_{T'}$ are also zero coordinates of $S_{T}$. An in-depth view of this partial order can be found in \cite{Postnikov}. We will later explore this partial order through the lense of bicolored tilings. For our purposes, it is sufficient to know that $\overline{S_T} = \overline{S_{T'}} \Longrightarrow S_T = S_{T'}$ (the converse is evidently true as well). Thus if we wanted to show that two positroid cells are the same, we may instead prove that their closures are the same (e.g. \Cref{localTransformationLemma}).
    
\end{definition}

\begin{example}
    We consider the tiling $T$ of type $(4,6)$ from the previous example. Let $I = 2456$. Then the matchings $m$ with $\partial m = I$ are
    \begin{center}
        \includegraphics[scale = 0.5]{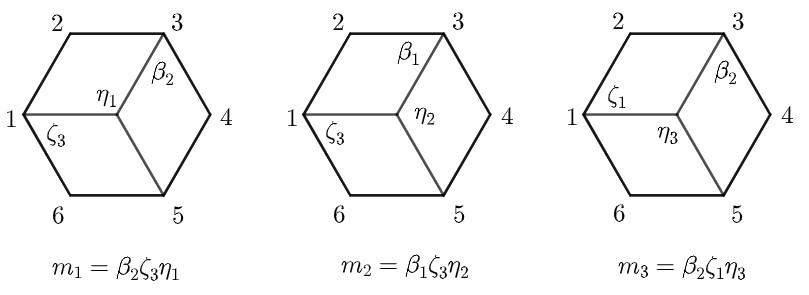}
    \end{center}
    Then $\Delta_{2456} = \beta_2 \zeta_3 \eta_1 + \beta_1 \zeta_3 \eta_2 + \beta_2 \zeta_1 \eta_3$.
    Doing this for every $k$-subset $I$ of $[n]$, we obtain
    \begin{align*}
        P_T &= (\Delta_{1234},\Delta_{1235},\Delta_{1236},\Delta_{1245}, \Delta_{1246},\Delta_{1256},\Delta_{1345},\\
        &\Delta_{1346},\Delta_{1356},\Delta_{1456},\Delta_{2345},\Delta_{2346},\Delta_{2356},\Delta_{2456},\Delta_{3456})
    \end{align*}
    with
    \begin{align*}
        \Delta_{1234} &= \delta_2 \varepsilon_3 \eta_1 \\
        \Delta_{1235} &= \gamma_2 \varepsilon_3 \eta_1 \\
        \Delta_{1236} &= \gamma_2 \delta_3 \eta_1 \\
        \Delta_{1245} &= \beta_1 \varepsilon_3 \eta_2 + \beta_2 \varepsilon_3 \eta_1 \\
        \Delta_{1246} &= \beta_1 \delta_2 \eta_3 + \beta_1 \delta_3 \eta_2 + \beta_2 \delta_3 \eta_1 \ \\
        \Delta_{1256} &= \beta_1 \gamma_2 + \eta_3 \\
        \Delta_{1345} &= \alpha_1 \varepsilon_3 \eta_2  \\
        \Delta_{1346} &= \alpha_1 \delta_2 \eta_3 + \alpha_1 \delta_3 \eta_2 \\
        \Delta_{1356} &= \alpha_1 \gamma_2 \eta_3 \\
        \Delta_{1456} &= \alpha_1 \beta_2 \eta_3 \\
        \Delta_{2345} &= \varepsilon_3 \zeta_1 \eta_2 \\
        \Delta_{2346} &= \delta_2 \zeta_1 \eta_3 + \delta_2 \zeta_3 \eta_1 + \delta_3 \zeta_1 \eta_2 \\
        \Delta_{2356} &= \gamma_2 \zeta_1 \eta_3 + \gamma_2 \zeta_3 \eta_1 \\
        \Delta_{2456} &= \beta_2 \zeta_3 \eta_1 + \beta_1 \zeta_3 \eta_2 + \beta_2 \zeta_1 \eta_3 \\
        \Delta_{3456} &= \alpha_1 \zeta_3 \eta_2
    \end{align*}
    And finally, $S_T = \{P_T \, | \, \alpha_1,\beta_1,\beta_2,\gamma_2,\delta_2,\delta_3,\varepsilon_3, \zeta_1,\zeta_3,\eta_1,\eta_2,\eta_3 > 0\}$.
\end{example}

\begin{definition}
    Let $G$ be a plabic graph. A \textit{perfect orientation} of $G$ is a choice of an orientation of its edges such that each black vertex has exactly one outgoing arrow, and each white vertex has exactly one incoming arrow. An \textit{almost perfect matching} of $G$ is a choice of edges of $G$ such that each internal vertex of is adjacent to exactly one edge in that subset of edges.
\end{definition}

We recall that if $T$ has $n$ boundary vertices, then $G = \Phi(T)$ has $n$ boundary vertices too, which are all of degree $1$ and incident to a single boundary edge.

\begin{remark}
    Let $T$ be a tiling, and $G = \Phi(T)$ the corresponding plabic graph (\Cref{stellarReplacementMap}). Any matching $m$ of $T$ gives rise to an almost perfect matching of $G$ and a perfect orientation of $G$. This works as follows:
    
    Any angle $\alpha$ of $T$ maps to a unique edge $e_{\alpha}$ of $G$: Let $m = \{\alpha_i \, | \, i \in \mathcal{I}\}$ be a matching of $T$. Let $\tilde{E}$ be the set of boundary edges in $G$ that are adjacent to white vertices in $G$ whose pre-image under $\Phi$ is in $\partial m$, i.e. they aren't matched in $m$. Then $\tilde{m} = \{ e_{\alpha_{i}} \, | \, i \in \mathcal{I} \} \cup \tilde{E}$ is an almost perfect matching of $G$.
    
    At the same time, an almost perfect matching results in a perfect orientation of $G$ by orienting all edges in $\tilde{m}$ from the black vertex to the white vertex, and all other edges that are not in $\tilde{m}$ the other way around.
\end{remark}

\begin{example}
    Consider the tiling $T$ from \Cref{runningExampleHex}, and the matching $m = \{\alpha_1, \gamma_2, \eta_3 \}$. For simplicity, we will label $a= \alpha_1$, $b = \gamma_2$, and $c= \eta_3$.
    \begin{center}
        \includegraphics[scale=0.5]{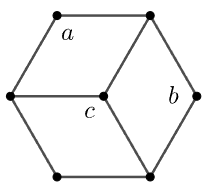}
    \end{center}
    Let $G = \Phi(T)$ be the corresponding plabic graph.
    \begin{center}
        \includegraphics[scale=0.5]{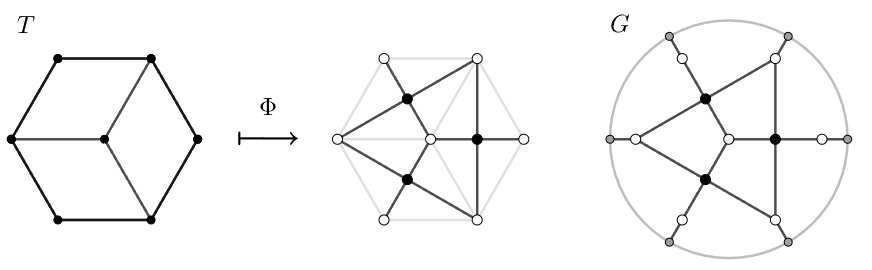}
    \end{center}
    Then the corresponding almost perfect matching of $G$ and perfect orientation of $G$ are
    \begin{center}
        \includegraphics[scale=0.5]{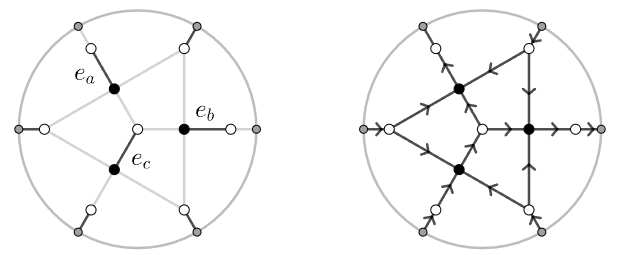}
    \end{center}
\end{example}

Using the allocation of an almost perfect matching to any tiling as described above, it follows that the parametrisation described in \Cref{tilingParametrisation} aligns with the parametrisation of positroid cells described in \cite[2.9-2.17]{Williams}, and consequently the parametrisation of Grassmann cells in \cite[Ch.11]{Postnikov}.

Next, we want to prove that equivalence, mutation, and reductions of a tiling preserve the corresponding positroid cell. To do so, we explore how local changes of a tiling change the positroid cell of the whole tiling.

\begin{remark}
    Let $T$ be a tiling of type $(k,n)$ of a disk $D_n$. Let $T'$ be a subtiling of $T$ under a disk $D_m \subset D_n$ as in \Cref{subtilingDefinition}, and $\tilde{T} = T \setminus T'$. We can extend the definition of \textit{matchings} \sout{described in Definition 4.2} to tilings of surfaces other than disks, such as $\tilde{T}$, i.e. we choose angles in $\tilde{T}$ satisfying the conditions (i)-(iii) in \Cref{matchingDefinition}.
    
    \begin{center}
        \includegraphics[scale=0.4]{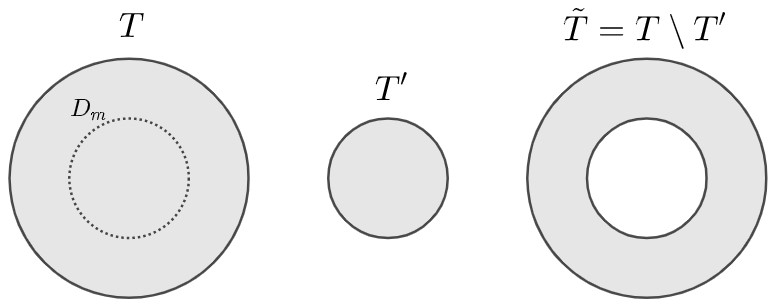}
    \end{center}
    
    We partition the boundary of a matching of $\tilde{T}$. If $\tilde{m}$ is a matching of $\tilde{T}$, we will denote $\partial_1{\tilde{m}}$ the set of boundary vertices of $\tilde{T}$ on the boundary of $T$ that \textit{are not} matched in $\tilde{m}$, and we will denote $\partial_2 \tilde{m}$ the boundary vertices of $\tilde{T}$ on the boundary of $T'$ that \textit{are} matched in $\tilde{m}$. Then, similarly, we define the $I$-th coordinate of $\tilde{T}$ as
    \[
        \tilde{\Delta_I} = \sum_{\substack{\tilde{m} \in \M (\tilde{T})\\ \partial_1 \tilde{m} = I}} \tilde{m}
    \]
    
    If $m$ is a matching of $T$, we  will write $m^{\circ}$ for the angles of $m$ inside $T'$, and $\tilde{m}$ for the angles of $m$ inside $\tilde{T}$, i.e. $m^{\circ} = m \restrict{T'}$ and $\tilde{m}= m \restrict{\tilde{T}}$. Note that in that case, $m = m^{\circ} \sqcup \tilde{m}$, and we have $\partial_1 \tilde{m} = \partial m$ and $\partial_2 \tilde{m} = \partial m^{\circ}$.
\end{remark}

\begin{example}
    Here is an example of a tiling $T$ of type $(8,13)$ and a matching $m \in \M(T)$.
    \begin{center}
        \includegraphics[scale=0.4]{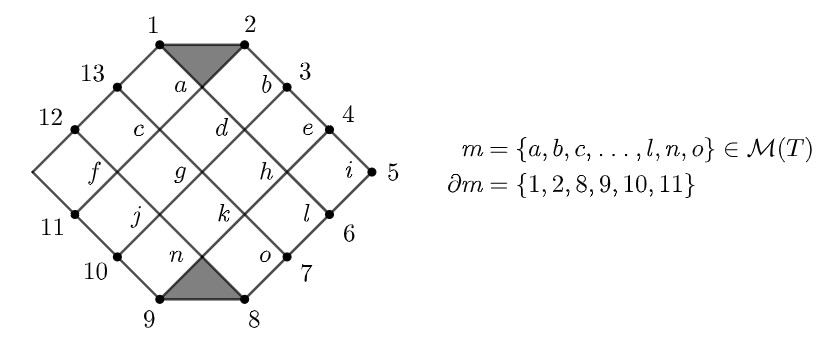}
    \end{center}
    Then we can "cut out" the subtiling $T'$ consisting of the four white tiles in the center of $T$. The remainder is the tiling $\tilde{T}$.
    \begin{center}
        \includegraphics[scale=0.4]{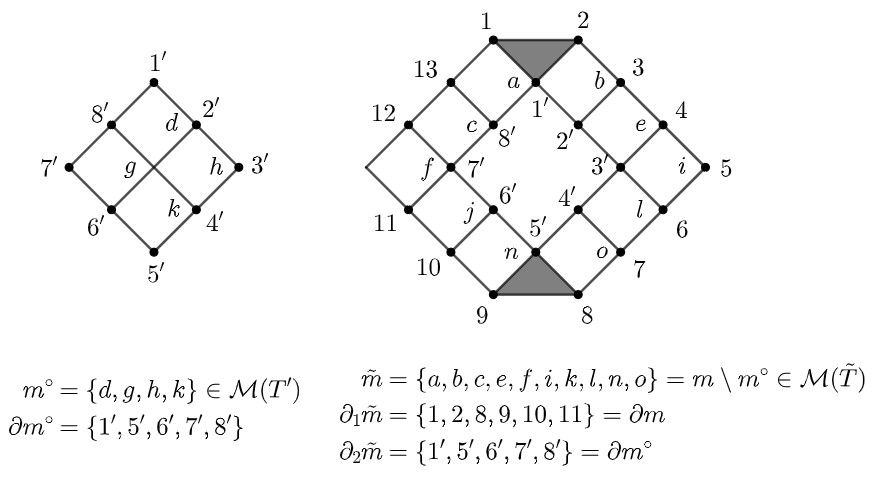}
    \end{center}
\end{example}

\begin{lemma} \label{localParametrisationFactorisation}
    Let $T$ be a tiling of type $(k,n)$. Let $T'$ be a subtiling of $T$, and $\tilde{T} = T \setminus T'$. Then for any $k$-subset $I$ of $[n]$, the $I$-th coordinate of the positroid cell $S_T$ is given by
    \[
        \Delta_I = \sum_{\substack{\tilde{m} \in \M(\tilde{T})\\ \partial_1 \tilde{m} = I}} \Delta_{\partial_2 \tilde{m}}^{\circ} \cdot \tilde{m}
    \]
    where $\Delta_J^{\circ}$ denotes the $J$-th coordinate of the positroid cell $S_{T'}$.
    \begin{proof}
        Let $m \in \M(T)$ be a matching of $T$ with $\partial m = I$. We write $m = m^{\circ} \sqcup \tilde{m}$, where $m^{\circ}$ denotes the angles of $m$ that are inside $T'$, and $\tilde{m}$ the angles outside $T'$ (that is, they are in $\tilde{T}$). Then $m^{\circ}$ is a matching of $T'$ and $\tilde{m}$ is a matching of $\tilde{T}$.
        
        Let $J = \partial m^{\circ}$. Then for any matching $m'$ of $T'$ with $\partial m' = J$, $m' \sqcup \tilde{m}$ is a matching of $T$ with $\partial (m' \sqcup \tilde{m}) = I$. In other words, all the matchings $m$ of $T$ with $\partial m = I$ are given as a combination of a matching $\tilde{m}$ of $\tilde{T}$ with $\partial_1 \tilde{m} = I$ and a matching $m^{\circ}$ of $T'$ with $\partial m^{\circ} = \partial_2 \tilde{m}$. By summing them as monomials, we obtain the $I$-th coordinate
        \[
            \Delta_I = \sum_{\substack{\tilde{m} \in \M (\tilde{T})\\ \partial_1 \tilde{m} = I}} \tilde{m} \cdot \left( \sum_{\substack{m^{\circ} \in \M (t)\\ \partial m^{\circ} = \partial_2 \tilde{m}}} m^{\circ} \right)
        \]
        where the second sum is the coordinate $\Delta_{\partial_2 \tilde{m}}^{\circ}$ of the positroid cell $S_{T'}$. Thus,
        \[
            \Delta_I = \sum_{\substack{\tilde{m} \in \M (\tilde{T})\\ \partial_1 \tilde{m} = I}} \tilde{m} \cdot \Delta_{\partial_2 \tilde{m}}^{\circ}
        \]
    \end{proof}
\end{lemma}

\begin{lemma} \label{localTransformationLemma}
    Let $T_1$ be a tiling with subtiling $A$. Let $T_2$ be the tiling obtained by replacing $A$ in $T_1$ with a new subtiling $B$ of same type as $A$ such that $\overline{S_{A}} \subset \overline{S_{B}}$. Then $\overline{S_{T_1}} \subset \overline{S_{T_2}}$.
    \begin{proof}
        Let $\alpha_1,\dots,\alpha_r$ and $\alpha_{r+1}, \dots, \alpha_R$ be the angles in $A$ and $\tilde{T} = T_1 \setminus A$, respectively. Let $\beta_1, \dots, \beta_s$ be the angles in $B$.
        
        Consider a point $x = P_{T_1}(x_1,\dots,x_r,x_{r+1}, \dots, x_R) \in \overline{S_{T_1}}$, with $x_i \geq 0$. We want to express $x$ as a point in $\overline{S_{T_2}}$, parametrised by $P_{T_2}$. We know that $x^{\circ} = P_A(x_1, \dots, x_r)$ is a point in $\overline{S_A}$. Since $\overline{S_A} \subset \overline{S_{B}}$, there are $y_1, \dots, y_s \geq 0$ such that $x = P_{B}(y_1,\dots,y_s) \in \overline{S_{B}}$.
        
        In other words, if $\Delta_J^A$ and $\Delta_J^B$ are the $J$-th coordinate in $P_A$ and $P_{B}$, respectively, then there is a $\lambda \in \R_{>0}$ such that
        \[
            \Delta_J^A(x_1,\dots,x_r) = \lambda \Delta_J^B(y_1,\dots,y_s) \,\,\,\,\,\,\,\,\,\,\,\,\,\,\,\,\,\,\,\, (\ast)
        \]
        Then for any $k$-subset $I$ of $[n]$, if $\Delta_I^{(1)}$ and $\Delta_I^{(2)}$ are the $I$-th coordinate in $P_{T_1}$ and $P_{T_2}$, we have by \Cref{localParametrisationFactorisation}
        \begin{align*}
            & \Delta_I^{(1)} (x_1,\dots,x_r,x_{r+1},\dots,x_R) & \\
            = \, & \sum_{\mathclap{\substack{\tilde{m} \in \M (\tilde{T})\\ \partial \tilde{m} = I}}} \, \Delta_{\partial_2 \tilde{m}}^A (x_1,\dots,x_r) \cdot \tilde{m}(x_{r+1},\dots,x_R) & \text{by \Cref{localParametrisationFactorisation}} \\
            = \, & \sum_{\mathclap{\substack{\tilde{m} \in \M (\tilde{T})\\ \partial \tilde{m} = I}}} \, \lambda \cdot \Delta_{\partial_2 \tilde{m}}^B (y_1,\dots,y_s) \cdot \tilde{m}(x_{r+1},\dots,x_R) & \text{using $(\ast)$} \\
            = \, & \lambda \cdot \sum_{\mathclap{\substack{\tilde{m} \in \M (\tilde{T})\\ \partial \tilde{m} = I}}} \,  \Delta_{\partial_2 \tilde{m}}^B (y_1,\dots,y_s) \cdot \tilde{m}(x_{r+1},\dots,x_R) & \\
            = \, & \lambda \Delta_I^{(2)} (y_1,\dots,y_r,x_{r+1},\dots,x_R) & \text{by \Cref{localParametrisationFactorisation}}
        \end{align*}
        Thus $x = (\Delta_I^{(2)}(y_1,\dots,y_s,x_{r+1},\dots,x_R))_{I \in \nSk} \in \overline{S_{T_2}}$. Hence $\overline{S_{T_1}} \subset \overline{S_{T_2}}$, which concludes the proof.
    \end{proof}
\end{lemma}

\begin{corollary} \label{localTransformationCorollary}
    Let $T_1$ be a tiling with subtiling $A$, and $T_2$ the tiling obtained by replacing $A$ in $T_1$ with a new subtiling $B$ with $\overline{S_A} = \overline{S_{B}}$. Then $\overline{S_{T_1}} = \overline{S_{T_2}}$. Equivalently, if $S_A = S_B$, then $S_{T_1} = S_{T_2}$.
\end{corollary}

\begin{proposition} \label{invariancesParametrisation}
    The positroid cell associated to $T$ is invariant under
    \begin{enumerate}[label = $(\roman*)$]
        \item mutation of a tiling.
        \item tiling equivalence.
        \item reductions of a tiling.
    \end{enumerate}
    \begin{proof}
    For all three parts of the proof we will pick a point in the open cell $S_{T_1}$ of a tiling $T_1$ and show that that point can be expressed as a point in the closed cell $\overline{S_{T_2}}$ where $T_2$ is the tiling obtained after transforming $T_1$ as described in $(i)$,$(ii)$, and $(iii)$, respectively, thus showing that $S_T \subset \overline{S_{T'}}$. It is crucial that we use the open positroid cell of $T'$, to avoid division by $0$ in some of the calculations. Since $S_T \subset \overline{S_{T'}} \Longrightarrow \overline{S_T} \subset \overline{S_{T'}}$, we still obtain the desired result.
    
    Moreover, by \Cref{localTransformationCorollary} it is sufficient to show that these transformations preserve the positroid cell locally.
    
    \begin{enumerate}[label = $(\roman*)$]
        \item
            Consider the two triangulations $T_1$ and $T_2$ of a quadrilateral and their parametrisations as described in \Cref{tilingParametrisation}.
            \begin{center}
                \includegraphics[scale=0.6]{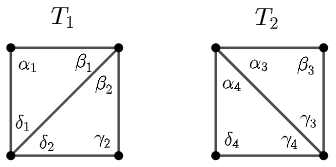}
            \end{center}
            \[
                P_{T_1} = (a_1 \beta_1, \alpha_1 \gamma_2, \alpha_1 \delta_2, \beta_1 \gamma_2, \beta_1 \delta_2 + \beta_2 \delta_1, \gamma_2 \delta_1)
            \]
            \[
                P_{T_2} = (\alpha_4 \beta_3, \alpha_3 \gamma_4 + \alpha_4 \gamma_3, \alpha_3 \delta_4, \beta_3 \gamma_4, \beta_3 \delta_4, \gamma_3 \delta_4)
            \]
            
            Let $x \in S_{T_1}$ given by the fixed parameters $a_1, b_1, b_2, c_2, d_1, d_2 > 0$, i.e.
            \[
                x = (a_1 b_1, a_1 c_2, a_1 d_2, b_1 c_2, b_1 d_2 + b_2 d_1, c_2 d_1)
            \]
            Then let $y \in \overline{S_{T_2}}$ be the point given by the fixed, non-negative parameters
            \begin{center}
                \begin{alignat*}{3}
                    \beta_3 &= a_1 && \delta_4 = \frac{b_1 d_2 + b_2 d_1}{a_1}\\
                    \alpha_4 &= b_2 && \alpha_3 = \frac{a_1 d_2}{\delta_4}\\
                    \gamma_4 &= \frac{b_1 c_2}{a_1} \,\,\,\,\,\,\, && \gamma_3= \frac{c_2 d_1}{\delta_4}
                \end{alignat*}
            \end{center}
            Then $y = (a_1 b_1, a_1 c_2, a_1 d_2, b_1 c_2, b_1 d_2 + d_2 d_1, c_2 d_1) = x$, thus $x \in \overline{S_{T_1}}$. Hence $S_{T_1} \subset \overline{S_{T_2}}$. Then we also have $S_{T_2} \subset \overline{S_{\mu_e(T_2)}} = \overline{S_{T_1}}$, and thus $S_{T_1} = S_{T_2}$.
        \item
            \begin{enumerate} [label = $(\alph*)$]
                \item (Hourglass equivalence)
                    Consider the tiling $T_1$ that is an empty $n$-gon, and $T_2$ that is obtained by adding an hourglass inside $T_1$.
                    \begin{center}
                        \includegraphics[scale = 0.4]{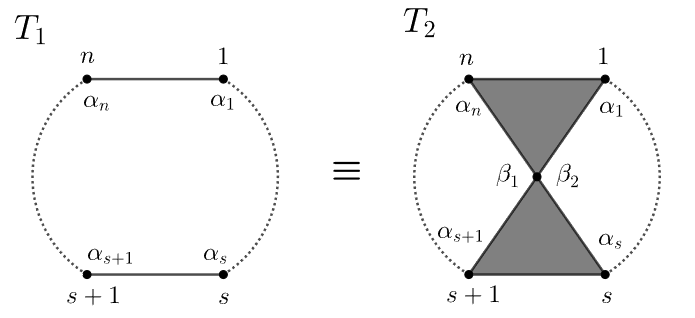}
                    \end{center}
                    
                    We call $\alpha_i$ the angle at the boundary vertex $i$ in both tilings, and $\beta_1$ and $\beta_2$ the angles around the internal vertex in $T_2$. The parametrisations of these tilings are
                    \[
                        P_{T_1} = (\alpha_n, \dots , \alpha_1)
                    \]
                    \[
                        P_{T_2} = (\alpha_n \beta_2, \dots, \alpha_{s+1} \beta_2, \alpha_{s} \beta_1, \dots, \alpha_1 \beta_1)
                    \]

                    \begin{enumerate}[label = $\cdot$]
                        \item
                            Let $x = (a_n,\dots,a_1) \in S_{T_1}$ with $a_i > 0$ for any $i = 1,\dots,n$.
                            We construct $y \in \overline{S_{T_2}}$, with parameters $\alpha_i = a_i \geq 0$ and $\beta_1=\beta_2=1$. Then
                            \[
                                y =  (a_n \cdot 1, \dots , a_1 \cdot 1) =x,
                            \]
                            thus $x \in \overline{S_{T_2}}$. Thus, $S_{T_1} \subset \overline{S_{T_2}}$.
                        \item
                            Let $y \in S_{T_2}$ with parameters $a_i$, $b_j>0$, that is
                            \[
                                y = (a_n b_2, \dots, a_{s+1} b_2, a_{s} b_1, \dots, a_1 b_1).
                            \]
                            We construct $x \in \overline{S_{T_1}}$ with parameters
                            \[
                                \alpha_i =
                                \begin{cases}
                                    a_ib_2, & \text{if } i = 1,\dots,s \\
                                    a_ib_1, & \text{if } i = s+1, \dots, n.
                                \end{cases}
                                \,\,\,\, \in \R_{\geq 0}
                            \]
                            Then $x = y$, and $y \in \overline{S_{T_1}}$. Thus, $S_{T_2} \subset \overline{S_{T_1}}$. 
                    \end{enumerate}
                    We conclude that $S_{T_1} = S_{T_2}$.
                \item
                    Consider the tiling $T$ and one of its boundary vertices $i$. Let $T'$ be the tiling obtained by decontracting at $i$, i.e.
                    \begin{center}
                        \includegraphics[scale = 0.35]{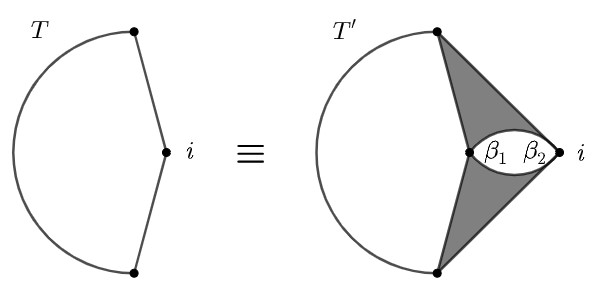}
                    \end{center}
                    We call $\alpha_1,\dots,\alpha_r$ the angles in the angles in $T$, with $\alpha_1,\dots,\alpha_s$ being the angles around vertex $i$ in $T$. We call $\beta_1, \beta_2$ the angles in the digon. Let $I$ be a $k$-subset of $[n]$. Let $\Delta_I$ be the $I$-th coordinate of $P_T$, and $\Delta_I'$ be the $I$-th coordinate of $P_{T'}$.
                    \begin{enumerate} [label = $\cdot$]
                        \item If $i \in I$, then for all matchings $m$ of $T$ with $\partial m = I$, $m' = m \sqcup \beta_1$ is a matching of $T'$, with $\partial m' = I$.
                        \item If $i \notin I$, then for all matchings $m$ of $T$ with $\partial m = I$, $m' = m \sqcup \beta_2$ is a matching of $T'$, with $\partial m' = I$.
                    \end{enumerate}
                    Thus by defining
                    \begin{center}
                        \begin{align*}
                            \lambda_I (\beta_1,\beta_2) = \begin{cases}
                                \beta_1 \text{ if } i \in I \\
                                \beta_2 \text{ if } i \notin I
                            \end{cases}
                        \end{align*}
                    \end{center}
                    we get that for matching $m$ of $T$, $m' = m \sqcup \lambda_{\partial m}$ is a matching of $T'$, and for any $I$, $\Delta_I' = \lambda_I \Delta_I$.
                    \begin{enumerate} [label = $\cdot$]
                        \item 
                            Let $x = P_T(a_1, \dots, a_r) \in S_T$. In other words, the $I$-th coordinate of $x$ is $\Delta_I(a_1, \dots, a_r)$. Then we construct $y \in \overline{S_{T'}}$ with parameters $\alpha_i = a_i$ and $\beta_i = 1$. Then the $I$-th coordinate of $y$ is
                            \[
                                \Delta_I'(a_1,\dots,a_r,1,1) = \lambda_I(1,1) \Delta_I(a_1,\dots,a_r) = \Delta_I(a_1,\dots,a_r)
                            \]
                            Thus $y = x$, and thus $x \in \overline{S_{T'}}$. Hence, $S_{T} \subset \overline{S_{T'}}$.
                        \item
                            Let $y = P_{T'}(a_1,\dots,a_r,b_1,b_2) \in S_{T'}$. In other words, the $I$-th coordinate of $y$ is $\Delta_I'(a_1, \dots, a_r,b_1,b_2)$. We recall that monomials appearing in any coordinate have the same length $p$ (i.e. the number of angles/variables in the monomial which equals the number of faces of the tiling). We denote $q = p^{-1}$. Then we construct $x = P_T(a_1 \mu, \dots, a_s \mu, a_{s+1}\nu, \dots, a_r \nu) \in \overline{S_{T}}$, where $\mu =  b_2 b_1^{-q(p-1)}$ and $\nu = b_1^q$. Then
                            \[
                                \Delta_I(a_1 \mu,\dots,a_s \mu,a_{s+1} \nu, \dots, a_r \nu) =  \sum_{\partial m = I} m(a_1 \mu,\dots,a_s \mu,a_{s+1} \nu, \dots, a_r \nu)
                            \]
                            Evaluating a monomial on the parameters equates to multiplying $p$ of the parameters (corresponding to the angles in the matching). We distinguish two cases
                            \begin{enumerate}[label = $-$]
                                \item $i \in I$. Then $\alpha_1,\dots,\alpha_s \notin m$. Thus the monomial is completely independent of those first $s$ parameters and is a product of $p$ of the remaining parameters. We can then write
                                \begin{center}
                                    \begin{align*}
                                       m(a_1 \mu,\dots,a_s \mu,a_{s+1} \nu, \dots, a_r \nu) &= \nu^p m(a_1,\dots,a_s,a_{s+1}, \dots, a_r) \\
                                       &= b_1 m(a_1, \dots, a_r) 
                                    \end{align*}
                                \end{center}
                                Thus
                                \begin{center}
                                    \begin{align*}
                                        \Delta_I(a_1 \mu,\dots,a_s \mu,a_{s+1} \nu, \dots, a_r \nu) &= b_1 \sum_{\partial m = I} m(a_1,\dots, a_r) \\
                                        &= \lambda_I(b_1,b_2) \Delta_I(a_1,\dots,a_r) \\
                                        &= \Delta_I'(a_1,\dots,a_r,b_1,b_2)
                                    \end{align*}
                                \end{center}
                                \item $i \notin I$. Then there is exactly one $j \in \{1,\dots,s\}$ such that $\alpha_j \in m$, which means exactly one copy of $\mu$ appears. The remaining $\{\alpha_1,\dots,\alpha_s\} \setminus \{\alpha_j\}$ do not appear in the monomial $m$, and instead $p-1$ of the angles $\alpha_{s+1}, \dots, \alpha_r$ do. Thus
                                \begin{center}
                                    \begin{align*}
                                       m(a_1 \mu,\dots,a_s \mu,a_{s+1} \nu, \dots, a_r \nu) &= \mu \nu^{p-1} m(a_1,\dots,a_s,a_{s+1}, \dots, a_r) \\
                                       &= b_2 m(a_1, \dots, a_r) 
                                    \end{align*}
                                \end{center}
                                Thus
                                \begin{center}
                                    \begin{align*}
                                        \Delta_I(a_1 \mu,\dots,a_s \mu,a_{s+1} \nu, \dots, a_r \nu) &= b_2 \sum_{\partial m = I} m(a_1,\dots, a_r) \\
                                        &= \lambda_I(b_1,b_2) \Delta_I(a_1,\dots,a_r) \\
                                        &= \Delta_I'(a_1,\dots,a_r,b_1,b_2)
                                    \end{align*}
                                \end{center}
                            \end{enumerate}
                            Thus for any $k$-subset $I$ of $[n]$, $\Delta'(a_1 \mu,\dots,a_s \mu,a_{s+1} \nu, \dots, a_r \nu) = \Delta_I'(a_1, \dots, a_r,b_1,b_2)$, thus $x=y$, and thus $y \in \overline{S_T}$. Hence $S_{T'} \subset \overline{S_T}$.
                    \end{enumerate}
                    We conclude that $S_T = S_{T'}$.
            \end{enumerate}
        \item
            Consider the tiling $T_1$ that is an $n$-gon with one $1$-edge $e$ (w.l.o.g at boundary vertex $n$). Let $T_2 = R_e(T_1)$, i.e. $T_2$ is an empty $n$-gon.
            \begin{center}
                \includegraphics[scale = 0.4]{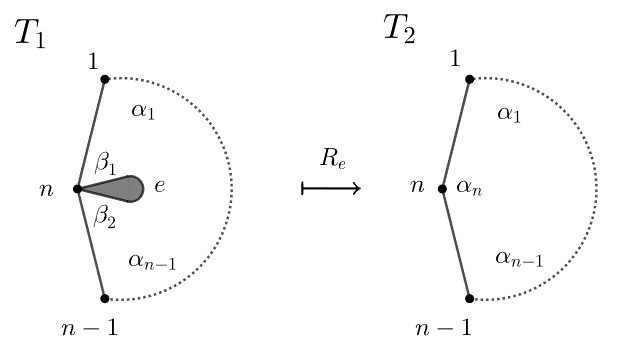}
            \end{center}
            We call $\alpha_i$ the angle at the boundary vertex $i$ in both tilings, except for the angles at the vertex $n$ in $T_1$ which we call $\beta_1$ and $\beta_2$. The parametrisations of these tilings are
            \[
                P_{T_1} = (\beta_1 + \beta_2, \alpha_{n-1}, \dots, \alpha_{1})
            \]
            \[
                P_{T_2} = (\alpha_n, \dots, \alpha_1)
            \]
            
            \begin{enumerate} [label = $\cdot$]
                \item
                    Let $x \in S_{T_1}$ with parameters $\alpha_i = a_i$, $\beta_j = b_j > 0$, that is
                    \[
                        x = (b_1 + b_2, a_{n-1}, \dots, a_{1})
                    \]
                    We construct $y \in \overline{S_{T_2}}$ with parameters $\alpha_n = b_1+b_2 \geq 0$, and $\alpha_i = a_i \geq 0$ for $i = 1,\dots, n-1$. Then $y = x$, and $x \in \overline{S_{T_2}}$. Thus, $S_{T_1} \subset \overline{S_{T_2}}$.
                \item
                    Let $y = (a_n, \dots, a_1) \in S_{T_2}$ with parameters $a_i >0$. We construct $x \in \overline{S_{T_1}}$ with parameters
                    \[
                        \begin{cases}
                            \alpha_i = a_i, & \text{for } i = 1,\dots,s-1 \\
                            \beta_1 = \beta_2 = \frac{1}{2}a_n
                        \end{cases}
                        \,\,\,\, \in \R_{\geq 0}
                    \]
                    Then $x=y$, and $y \in \overline{S_{T_1}}$. Thus $S_{T_2} \subset \overline{S_{T_1}}$.
            \end{enumerate}
            Hence, $S_{T_1} = S_{T_2}$.
    \end{enumerate}
    \end{proof}
\end{proposition}

\begin{theorem}
    Reduced (bicolored) tilings of type $(k,n)$ up to tiling equivalence are in bijection with positroid cells of the totally non-negative Grassmannian $\Grkn^{\tnn}$.
    \begin{proof}
        This follows from the fact that positroid cells are in bijection with Postnikov diagrams $\cite[14.2,14.7]{Postnikov}$ up to geometric exchange, which are in bijection with reduced tilings up to tiling equivalence.
    \end{proof}
\end{theorem}

\begin{proposition} \label{edgeIsDim}
    Let $E$ be the number of edges in a reduced tiling $T$. Then $\dim S_T = E - 1$.
    \begin{proof}
        Edges of $T$ map to faces of the plabic graph $G=\Phi(T)$. If $\mathcal{F}$ is the number of faces in $G$, then $\dim S_T = \mathcal{F} - 1$ by \cite[12.7]{Postnikov}, and thus $\dim S_T = E - 1$.
    \end{proof}
\end{proposition}

\section{Degenerations of tilings}

We recall from \Cref{tilingParametrisation} that the closure of a positroid cell is given by
\[
    \overline{S_T} = \{(\Delta_I) \mid \alpha \geq 0 \,\,\, \forall \alpha \in A \}
\]
We can describe this order in terms of tilings by defining the degeneration of tilings.

\begin{definition}
    We define a partial order on $\tilingSet$ by
    \[
        T < T' \Longleftrightarrow S_T \subset \overline{S_{T'}} \Longleftrightarrow \overline{S_T} \subset \overline{S_{T'}}
    \]
\end{definition}

Degenerations of a tiling happen with respect to angles of that tiling. In order to get consistent results, we distinguish between two types of angles as follows.

\begin{definition}
    An angle $\alpha$ of a tiling $T$ is said to be \textit{essential} if for any matching $m \in \M (T)$, we have $\alpha \in m$. Similarly, $\alpha$ is said to be \textit{non-essential} if there is a matching $m \in \M (T)$ such that $\alpha \notin m$.
\end{definition}

\begin{definition}
    Let $T$ be a tiling, and $\alpha \in A$ a non-essential angle of $T$. Let $v=v(\alpha)$ be the vertex at $\alpha$, and $f = f(\alpha)$ the face in which $\alpha$ lies. Let $e_1$ and $e_2$ be the two edges adjacent to $\alpha$, and let $v_1$ and $v_2$ be the second endpoints of $e_1$ and $e_2$, respectively. Let $T'$ be the tiling obtained by constructing a black triangle with endpoints $v$,$v_1$, and $v_2$ inside $f$, such that the edges $e_1,e_2$ of $T$ merge with the black triangle.
    \begin{center}
        \includegraphics[scale=0.5]{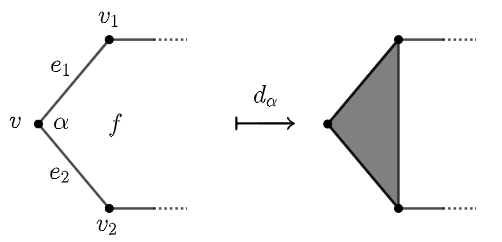}
    \end{center}
    Then $T'$ is called the \textit{degeneration of $T$ with respect to $\alpha$}, and is denoted $d_{\alpha}{T}$.
\end{definition}

\begin{proposition}
    The type of a tiling $T$ is invariant under degeneration.
    \begin{proof}
        Let $T$ be a tiling of type $(k,n)$ and $\alpha$ an angle of $T$. Let $T'= d_{\alpha}T$. Let $m$ be a matching of $T$ such that $\alpha \notin m$. Then $m$ is a matching of $T'$. Thus $T'$ has rank $\partial m = k$. Since the number of boundary vertices are not changed when degenerating, the type of $T'$ is $(k,n)$.
    \end{proof}
\end{proposition}

If $T$ is a tiling with diagram $S(T)$, then any intersection between two strands in $S(T)$ determines an angle. This follows from the definition of the Scott map (see \Cref{ScottMap}). If two strands $\gamma_i$ and $\gamma_j$ intersect and determine the angle $\alpha$, we denote $\alpha = \gamma_i \wedge \gamma_j$. Since strands may intersect more than once, we choose $\alpha$ to be the last intersection between $\gamma_i$ and $\gamma_j$ when following the orientation of $\gamma_i$.

\begin{proposition}
    Let $T$ be a tiling of permutation $\pi$, and $\gamma_i$, $\gamma_j$ be two distinct, intersecting strands of $S(T)$. Let $\alpha = \gamma_i \wedge \gamma_j$ be a non-essential angle of $T$. Let $\pi'$ be the permutation of $T' = d_{\alpha}T$. Then $\pi' = (\pi(i) \, \pi(j)) \circ \pi$.
    \begin{proof}
        This result is immediate if we observe how degenerations affect the diagram locally from $T$ to $T'$.
        \begin{center}
            \includegraphics[scale=0.5]{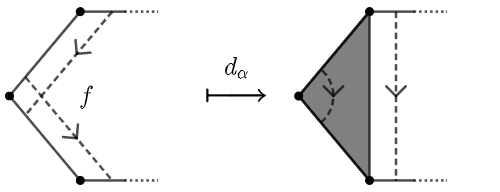}
        \end{center}
    \end{proof}
\end{proposition}

\begin{proposition} \label{degenerationParametrisation}
    Let $T$ be a tiling with angles $\alpha_1, \dots, \alpha_m$. Let $P_T = P_T(\alpha_1, \dots, \alpha_m)$. Let $T' = d_{\alpha_i}T$ for some non-essential angle $\alpha_i \in A$. Then the parametrisation of $S_{T'}$ is
    \[
        P_{T'} = P_{T'}(\alpha_1, \dots, \alpha_{i-1},\alpha_{i+1},\dots, \alpha_m) = P_T \restrict{\alpha_i = 0} = P_T(\alpha_1, \dots, \alpha_{i-1}, 0 , \alpha_{i+1}, \dots, \alpha_m)
    \]
    \begin{proof}
        Let $m$ be a matching of $T$ with $\alpha_i \notin m$. Then $m$ is also a matching of $T'$. Conversely, if $m'$ is a matching of $T'$, then $m'$ is a matching of $T$ as well, with $\alpha_i \notin m'$. In other words, the matchings of $T'$ are exactly the matchings of $T$ that do not contain $\alpha_i$. Thus, if $\Delta_I$ and $\Delta_I'$ denote the $I$-th coordinate of $T$ and $T'$, respectively, we get
        \[
            \Delta_I' = \sum_{\substack{m \in \M (T')\\ \partial m = I}} m
            = \sum_{\substack{m \in \M (T)\\ \partial m = I \\ \alpha_i \notin m}} m
            = \sum_{\substack{m \in \M (T)\\ \partial m = I \\}} m \restrict{\alpha_i = 0}
            = \left(\sum_{\substack{m \in \M (T)\\ \partial m = I \\}} m \right) \restrict{\alpha_i = 0}
            = \Delta_I \restrict{\alpha_i = 0}
        \]
        Hence
        \[
            P_{T'} = (\Delta_I') = (\Delta_I \restrict{\alpha_i = 0}) = P_T \restrict{\alpha_i=0}
        \]
        which concludes the proof.
    \end{proof}
    
\end{proposition}

\begin{corollary}
    If $T' = d_{\alpha} T$, then $T' < T$.
    \begin{proof}
        This follows directly from \Cref{degenerationParametrisation}, as $\overline{S_{T'}} \subset \overline{S_{T}}$.
    \end{proof}
\end{corollary}

The number of edges is reduced by exactly $1$ after degenerating a tiling, as two edges are merged together into one by adding a black triangle. From \Cref{edgeIsDim} the following immediately follows.

\begin{corollary}
    Let $T' = d_{\alpha}(T)$ be the degeneration of a tiling $T$ with respect to the angle $\alpha$. Then $\dim S_{T'} \leq \dim S_{T}-1$.
\end{corollary}

\begin{remark}
    The reason why we don't have an equality $\dim S_{T'} = \dim S_{T} - 1$ is that that the equality $\dim S_T = E - 1$ is only true if $T$ is reduced. After degenerating, the resulting tiling is not necessarily reduced.
\end{remark}

\begin{example}
    The following reduced tiling $T$ of type $(3,6)$ and of dimension $7$ can be degenerated at $\alpha$. The resulting tiling $T'$ is not reduced. After reducing $T'$ to a tiling $T''$, we see that the dimension of the corresponding positroid cell is $\dim S_{T''} = 5$.
    \begin{center}
        \includegraphics[scale=0.6]{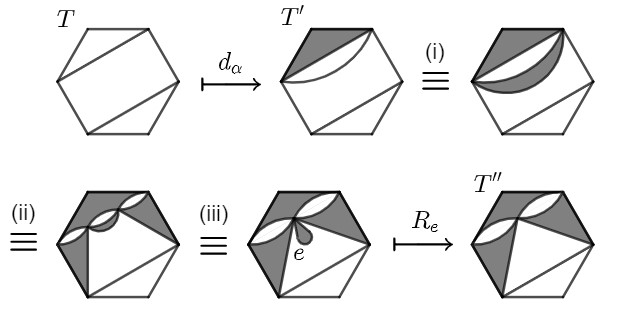}
    \end{center}
    The steps applied to the degenerated tiling $T'$ are as follows
    \begin{enumerate} [label = (\roman*)]
        \item Any simple edge is also a black digon.
        \item We decontract two white digons. It may be easier to see the transformation from right to left, by contracting the two digons that are adjacent to the boundary.
        \item We contract the central digon. This transforms the black digon into a $1$-gon that can be reduced to arrive at $T''$.
    \end{enumerate}
    
\end{example}

We summarise the main results of this section.

\begin{theorem}
    Let $T$ be a tiling of type $(k,n)$ with permutation $\pi$, and let $\alpha \in A$ such that $\exists$ distinct $i,j \in [n]$ with $\gamma_i \wedge \gamma_j = \alpha$. Let $T' \vcentcolon = d_{\alpha}(T)$ be the degeneration of $T$ at $\alpha$, and let $P=P(\alpha)_{\alpha \in A}$ be the parametrisation of $T$. Then
    \begin{enumerate}[label = $\bullet$]
        \item $T'$ is of type $(k,n)$.
        \item $T'$ has decorated permutation $\pi' = (\pi(i) \,\, \pi(j)) \circ \pi$.
        \item $T'$ parametrises the positroid cell $S_{\pi'}$ by $P \restrict{\alpha = 0}$.
        \item $T < T'$ and $\dim S_T < \dim S_{T'}$.
    \end{enumerate}
\end{theorem}

\newpage

\bibliographystyle{plain}
\bibliography{main}

\noindent
\begin{footnotesize}
\sc School of Mathematics, University of Leeds, Leeds LS2 9JT, UK\\
\textit{E-mail address:} \href{mailto:mm17jcdr@leeds.ac.uk}{\nolinkurl{mm17jcdr@leeds.ac.uk}},
\href{joelcosta94i@gmail.com}{\nolinkurl{joelcosta94i@gmail.com}}
\end{footnotesize} 

\end{document}